# Cubic Observers for a Class of Nonlinear Systems with Input and output delays

Mohammad Mahdi Share Pasand[1]

[1] Department of Electrical and Electronics, Faculty of Engineering & Technology, Standard Research Institute, Alborz, Iran, PO Box 31585-163, Email: sharepasand@standard.ac.ir

*Abstract.* This paper utilizes the recently proposed cubic observer to estimate the state of a class of nonlinear systems. The cubic observer is proposed as an alternative to linear observers for improved convergence rate and robustness. It is shown that the proposed cubic observer can enhance the estimation compared to linear observers used for the same process. The proposed observer can be viewed as a generalization to the linear observer. To further enhance the cubic observer, a composite observer is proposed which combines the high gain observer design principle with the cubic observer structure. An example is included to show the efficacy of the proposed observer.

## 1. Introduction

Observers have been the focus of decades of research in feedback control theory, condition monitoring and fault detection. Numerous aspects of observers have been investigated in the literature making the field a mature and sophisticated major. Different observer schemes have been developed to estimate the state of different linear systems including delayed [1-3], unknown input [4-6], output constrained [7], descriptor [8] and time varying [9, 10] linear systems. Design of observers for nonlinear systems has been also considered in the literature. Sliding mode [11] and high gain [12] observers are among nonlinear observer schemes used for state estimation of nonlinear systems.

Recently, a nonlinear observer, known as the cubic observer is proposed for state estimation of linear and systems [13]. The cubic observer contains a cubic term constructed by the product of a quadratic form and a linear multiplier of the output estimation error. This observer is distinguished from other nonlinear observers by its fast convergence especially when the initial estimation error is large [1]. The cubic observer is shown to yield improved performance while preserving desired robustness properties when used for linear systems [13]. The cubic term makes the Lyapunov function derivative of the estimation error dynamics of this observer to lie beneath that of a linear estimation error dynamics. Unlike the high gain observers [14], the estimation error dynamics of the cubic observer is nonlinear. Although its actual advantage is yet to be known in practice, the cubic observer is suggested as an alternative when other observer structures does not result satisfactory performance.

In this paper, we show that the cubic observer can be used for state estimation of Lipschitz and one-sided Lipschitz nonlinear systems as well. Lipschitz nonlinear systems have been the focus of several research studies [16] and proved to be capable of modeling many practical systems [17-23]. One-sided Lipschitz systems [24-29] form a wider class of nonlinear systems and include Lipschitz systems as a special sub-class. In the present paper, we assume a general nonlinear function consisting input, state and output delays. However, it is assumed that the non-Lipschitz part of the nonlinear function enters the process dynamics through a gain matrix which can be cast out by unknown input observer design techniques [4]. The Lipschitz part of the nonlinear function is replaced by a one-sided Lipschitz function in the second theorem. The one-sided Lipschitz condition covers a wider class of nonlinear systems and therefore, the second theorem generalizes the first. It should be noted that the expressions are however different.

This paper is organized in three sections. The second section provides the main result while the third and fourth sections are dedicated to simulation and conclusions.

## 2. Main result

### 2.1. Convergence Criteria

Assume a nonlinear system with multiple input delays, described by the following set of state space equations:

$$\dot{x}(t) = Ax(t) + f\left(x(t), u(t-\delta_1), .., u(t-\delta_{m_u}), y(t-\tau_1), ..., y(t-\tau_{m_y})\right)$$
$$y(t) = Cx(t) \tag{1}$$

Vector $x(t) \in R^n$ represents the state variable, $y(t) \in R^{n_y}$ is the measured output and matrices $A_{n \times n}$ and $C_{n_y \times n}$ are the state and output matrices respectively. Scalars $m_u, m_y$ represent the number of different time delays $\delta_1, ..., \delta_{m_u}, \tau_1, ..., \tau_{m_y}$ affecting the input and output vectors present in (1). By assigning zero values to the delays, the system becomes a delay-free one to which the results of this paper apply as well. Vector function $f(.)$ is assumed to have the following form:

$$f(x(t), U_\delta(t), Y_\tau(t)) = f_u(U_\delta(t), Y_\tau(t)) + Df_g(x(t), U_\delta(t), Y_\tau(t)) + f_L(x(t), U_\delta(t), Y_\tau(t)) \tag{2}$$

In which:

$$U_\delta(t) = \begin{bmatrix} u(t-\delta_1) \\ \vdots \\ u(t-\delta_{m_u}) \end{bmatrix}, \quad Y_\tau(t) = \begin{bmatrix} y(t-\tau_1) \\ \vdots \\ y(t-\tau_{m_y}) \end{bmatrix} ; \delta_1, ..., \delta_{m_u}, \tau_1, ..., \tau_{m_y} \geq 0 \tag{3}$$

Function $f(.)$ is formed by adding three nonlinear functions. A general nonlinear function of delayed inputs and outputs represented by $f_u(.)$, another general nonlinear function which depends not only on the delayed inputs and outputs but also on the actual state vector of the process represented by $f_g(.)$ and assumed to enter the process dynamics through an input distribution matrix $D$ which has to fulfill a specific property (See (10)) and a third nonlinear function $f_L$ which is assumed to be Lipschitz in $x(t)$ which means:

$$\exists \gamma > 0; \; \|f_L(x_1, U_\delta, Y_\tau) - f_L(x_2, U_\delta, Y_\tau)\| \leq \gamma \|x_1 - x_2\| ; \; \forall x_1, x_2, U_\delta, Y_\tau \tag{4}$$

Where the positive scalar $\gamma$ is the Lipschitz constant and $\|.\|$ represents the second norm of a vector. It is important to note that (4) assumes globally Lipschitz functions leading to global stability of the estimation error dynamics. Locally Lipschitz functions may be assumed, when local stability is considered. The presumed nonlinear function in the process dynamics (1) is a generalization to Lipschitz nonlinear systems previously considered in [9-12] and slightly enhances those in [13-15]. In order to estimate the state of (1) from the measured output, we propose a cubic observer, described by the following:

$$\dot{w}(t) = Gw(t) + Jy(t) + (I - EC)f_u(U_\delta(t), Y_\tau(t)) + (I - EC)f_L(\hat{x}(t), U_\delta(t), Y_\tau(t)) - e_c^T(t)C^T\theta Ce_c(t)NCe_c(t)$$
$$\hat{x}(t) = w + Ey(t) \tag{5}$$

Where the vector $w(t) \in R^{n_w}$ is an intermediate state estimation vector and $e_c(t) = x(t) - \hat{x}(t)$ is the estimation error of the cubic observer. Matrices $J, E, N$ and $\theta$ are observer parameters with appropriate dimensions to be determined in the design procedure. Vector $\hat{x}(t) \in R^n$ is the estimated state vector to be retrieved from the information in $w(t)$ and $y(t)$. In (5), a cubic term is added to the linear observer.

The output estimation error $Ce_c(t)$ is used to construct the cubic term. Observer (5) generalizes the previously studied observers for partial state estimation of nonlinear systems studied for instance in [15] chapter 6. The following theorem states convergence criteria for the proposed cubic observer.

***Theorem 1.*** The estimation error dynamics of observer (5) for system (1) is globally stable if (6)-(10) hold:

$$\nexists\, v \neq 0; \quad Gv + v^T C^T \theta C v N C v = 0 \tag{6}$$

$$\exists\, P = P^T > 0, \beta > 0; \begin{cases} \begin{bmatrix} PG + G^T P + \gamma^2 \beta I & P(I - EC) \\ (I - EC)^T P & -\beta I \end{bmatrix} < 0 \tag{7} \\ \\ PNC + C^T N^T P < 0 \tag{8} \end{cases}$$

$$\begin{cases} (I - EC)A - JC - G(I - EC) = 0 \tag{9} \\ (I - EC)D = 0 \tag{10} \end{cases}$$

In which $\theta = \theta^T \geq 0$.

***Proof:*** Before proceeding to the detailed proof, note that (6) guarantees that the nonlinear (cubic) estimation error dynamics does not possess any equilibrium point except the origin (i.e. $e_c = 0$). Conditions (7) and (8) are Lyapunov stability criteria for the origin as an equilibrium for the estimation error dynamics. Conditions (9) and (10) guarantee that the estimation error dynamics is not affected by the inputs, state and nonlinearities. Write the dynamics of $e_c(t)$ as;

$$\dot{e}_c(t) = \dot{x}(t) - \dot{\hat{x}}(t) = (I - EC)\dot{x}(t) - \dot{w}(t)$$

Substitute for $\dot{w}(t)$ from (5) and $\dot{x}(t)$ from (1). Use the system output equation in (1) and the estimation error definition to substitute for $w(t)$ and $y(t)$ respectively. The resulting equation reads;

$$\dot{e}_c(t) = Ge_c(t) + \big((I - EC)A - JC - G(I - EC)\big)x(t) + \cdots$$

$$(I - EC)\big(f_L(x(t), U_\delta, Y_\tau) - f_L(\hat{x}(t), U_\delta, Y_\tau)\big) + (I - EC)Df_g(x(t), U_\delta, Y_\tau) + e_c^T(t)C^T \theta C e_c(t) N C e_c(t)$$

In order to eliminate the effect of the state, delayed inputs and outputs and nonlinearities on the estimation error dynamics, it is sufficient that (8) hold. The estimation error dynamics can be simplified as;

$$\dot{e}_c(t) = Ge_c(t) + (I - EC)\Delta f_L + e_c^T(t) C^T \theta C e_c(t) N C e_c(t) \tag{11}$$

Where;

$$\Delta f_L = f_L(x(t), U_\delta, Y_\tau) - f_L(\hat{x}(t), U_\delta, Y_\tau) \tag{12}$$

Notice that the origin (i.e. $e_c = 0$) is a trivial equilibrium to (9) due to the Lipschitz property of $f_L(.)$ stated in (4). If condition (6) holds, then (11) have no equilibrium point except the origin. Therefore, stability of the origin implies global stability for (11). Now define the following Lyapunov function;

$$V = e_c^T P e_c \tag{13}$$

Taking derivative of (11);

$$\dot{V} = e_c^T(t)(PG + G^TP)e_c(t) + e_c^T(t)C^T\theta C e_c(t).e_c^T(t)(PNC + C^TN^TP)e_c(t) + e_c^T(t)P(I - EC)\Delta f_L$$
$$+ \Delta f_L^T(I - EC)^T P e_c(t)$$

Due to (4), for any arbitrary scalar $\beta > 0$, we have;

$$\beta \Delta f_L^T \Delta f_L \leq \gamma^2 \beta e_c^T(t) e_c(t)$$

Therefore;

$$\dot{V} \leq e_c^T(t)(PG + G^TP)e_c(t) + e_c^T(t)C^T\theta C e_c(t).e_c^T(t)(PNC + C^TN^TP)e_c(t) + e_c^T(t)P(I - EC)\Delta f_L$$
$$+ \Delta f_L^T(I - EC)^T P e_c(t) + \gamma^2 \beta e_c^T(t) e_c(t) - \beta \Delta f_L^T \Delta f_L$$

Rearrange the above inequality into matrix form to obtain;

$$\dot{V} \leq [e_c^T(t) \quad \Delta f_L^T]\begin{bmatrix} PG + G^TP + \gamma^2\beta I & P(I - EC) \\ (I - EC)^T P & -\beta I \end{bmatrix}\begin{bmatrix} e_c(t) \\ \Delta f_L \end{bmatrix} + e_c^T(t)C^T\theta C e_c(t) e_c^T(t)(PNC + C^TN^TP)e_c(t)$$

The term $e_c^T(t)C^T\theta C e_c(t)$ is a non-negative scalar since $\theta$ is assumed to be a positive semi-definite matrix. Therefore, in order for the Lyapunov function derivative to be negative, it is sufficient that the matrices in the core of the quadratic terms are negative definite. ∎

In order to widen the family of the nonlinear systems under discussion, one may replace Lipschitz functions by one sided Lipschitz functions. One-sided Lipschitz condition is less restrictive than Lipschitz condition and captures many nonlinear systems of interest [23]. Recall the following properties.

$$\exists \rho \in R \; ; \; \langle f_L(x_1, U_\delta, Y_\tau) - f_L(x_2, U_\delta, Y_\tau), x_1 - x_2 \rangle \leq \rho \|x_1 - x_2\|^2 \; ; \; \forall x_1, x_2, U_\delta, Y_\tau \tag{14}$$

$$\|f_L(x_1, U_\delta, Y_\tau) - f_L(x_2, U_\delta, Y_\tau)\|^2 \leq a\|x_1 - x_2\|^2 + b\langle f_L(x_1, U_\delta, Y_\tau) - f_L(x_2, U_\delta, Y_\tau), x_1 - x_2 \rangle \tag{15}$$

In which $\rho$ (the one sided Lipschitz constant), $a$ and $b$ are constants. Condition (15) is known as quadratic inner bounded condition [23, 26].

***Theorem 2.*** Assume that $f_L(.)$ fulfills (14), (15) instead of (4). The estimation error dynamics of observer (5) for system (1) is globally stable if (16)-(20) hold:

$$\nexists v \neq 0; \; Gv + v^T C^T \theta C v N C v = 0 \tag{16}$$

$$\exists P = P^T > 0, \mu_1, \mu_2 > 0; \begin{cases} \begin{bmatrix} PG + G^TP + 2(\mu_1\rho + \mu_2 a)I & (\mu_2 b - \mu_1)P(I - EC) \\ (\mu_2 b - \mu_1)(I - EC)^T P & -2\mu_1 I \end{bmatrix} < 0 & (17) \\ \\ PNC + C^T N^T P < 0 & (18) \end{cases}$$

$$\begin{cases} (I - EC)A - JC - G(I - EC) = 0 & (19) \\ (I - EC)D = 0 & (20) \end{cases}$$

In which $\theta = \theta^T \geq 0$.

***Proof.*** The proof is similar to that of *Theorem 1* with (14) and (15) replacing (4). ∎

In practical applications, it is inevitable to consider modeling uncertainties. For this purpose, assume the function $f_L(.)$ is uncertain belonging to a family $f_L(\varepsilon)$ fulfilling (4) or (14), (15) for all $\varepsilon$. Therefore, *Theorem 1* holds for all $\varepsilon$. Furthermore, if (14), (15) hold for all uncertainties captured by $f_L(\varepsilon)$, then *Theorem 2* holds for all $\varepsilon$. This means that *Theorem 1, 2* are valid in presence of any uncertainties in $f_L(.)$, as long as the corresponding Lipschitz conditions hold in presence of those uncertainties. Furthermore, we may combine the high gain observer design technique [12] to compensate for the effects of model mismatch.

## 2.2. Design Aspects

The following theorem from [13] provides a solution to (6) and (8) or (16) and (18). The proof is very similar to that presented in [13] with small differences in nomenclature. Therefore, it is omitted here.

***Theorem 3.*** [13] For a given $C^T \theta C > 0$, (21) simultaneously solves (6) and (8) as well as (16) and (18).

$$N = -\alpha P^{-1} C^T \theta \qquad (21)$$

Scalar $\alpha > 0$ is arbitrary. ∎

Conditions (10) and (20) can be fulfilled if and only if $rank(D) = rank(CD)$ [15]. Conditions (9) and (19) can be fulfilled by design procedures for unknown input observers including those in [1, 6] and references therein. Conditions (7) and (17) can be numerically examined to find $\beta$ or $\mu_1, \mu_2$.

## 3. Simulation example

Consider following system which is the form of (1).

$$f_u = \begin{bmatrix} u(t-1) \\ u(t) \end{bmatrix}, \quad f_L = \begin{bmatrix} x_1(t)\cos(u(t)) \\ \sin(x_2(t)) \end{bmatrix}, \quad f_g = x_2(t)x_1(t)$$

$$A = \begin{bmatrix} -2 & -10 \\ 0 & -1 \end{bmatrix}, \quad C = \begin{bmatrix} 1 & 0 \end{bmatrix}, \quad D = \begin{bmatrix} -1 \\ 1 \end{bmatrix}$$

To fulfill (9) and (10);

$$E = \begin{bmatrix} 1 \\ -1 \end{bmatrix}, \quad J = \begin{bmatrix} 0 \\ 9 \end{bmatrix}, \quad G = \begin{bmatrix} -10 & 0 \\ 1 & -11 \end{bmatrix}$$

Assume $\beta = 100$, $\gamma = 1$. A solution to (7) can be obtained by solving a convex optimization problem. Here, we used *CVX* toolbox of [21] defining the maximum eigenvalue of the left hand side matrix in (7) as the objective function to be minimized. The resulting matrix is;

$$P = \begin{bmatrix} 59.0535 & 1.7898 \\ 1.7898 & 17.8858 \end{bmatrix}$$

Using (21) with $\theta = I_2, \alpha = 1$, a solution to (8) is given by: $N = [-.017 \quad .0017]^T$. Define the cumulative estimation error as;

$$J_o(t) = \int_0^t \left(x_1(\tau) - \hat{x}_1(\tau)\right)^2 + \left(x_2(\tau) - \hat{x}_2(\tau)\right)^2 d\tau$$

The simulation is run with $\hat{x}_1(0) = \hat{x}_2(0) = -5$ and $x_1(0) = x_2(0) = 0$ for an observer with cubic estimation error dynamics ($\gamma = 1$) and another with linear estimation error dynamics (i.e. $\gamma = 0$).

In the first case, shown in **Fig.1** the plant model matches those of the observers. In the second case, the plant model is replaced by the following model, while the observers use the original model.

$$\tilde{A} = \begin{bmatrix} -.9 & -8.9 \\ 1.1 & .1 \end{bmatrix}, \qquad \tilde{f}_u = \begin{bmatrix} u(t-2) \\ u(t) \end{bmatrix}$$

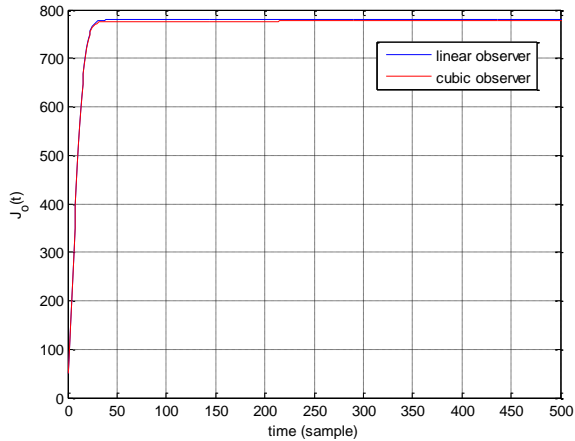
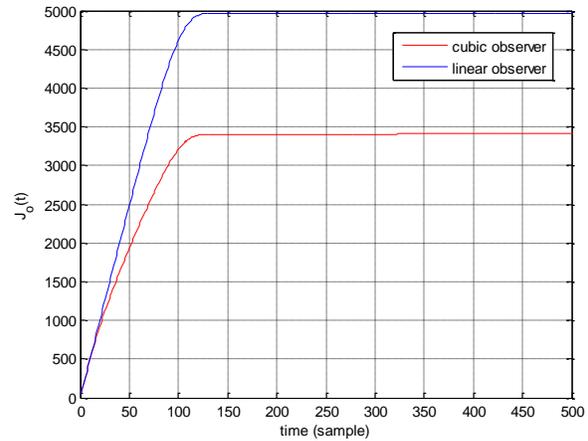

**Figure 1.** Cumulative estimation errors of cubic and linear estimation error dynamics; nominal model

**Figure 2.** Cumulative estimation errors of cubic and linear estimation error dynamics; uncertain model

In the uncertain model case, the cubic observer significantly outperforms a linear observer with the same observer parameters. In the nominal model case, the performance advantage is negligible.

## Conclusion

In this paper the use of a recently proposed nonlinear observer (the cubic observer) in state estimation of a class of nonlinear systems is examined. The considered class of nonlinear systems includes one-sided Lipschitz systems as a special sub-family. It is shown that if some general convergence criteria are fulfilled, the cubic observer can estimate the state of the given plant. In the simulation example, the cubic observer yielded improved performance when model uncertainties are assumed. Yet the actual usefulness of cubic observers requires future investigation, this paper suggests these observers to be considered in applications where linear observers fail to perform satisfactorily.